\documentclass[a4paper,12pt,twoside]{article}
\usepackage{amsmath,amsthm,amsfonts,amssymb}
\usepackage[latin1]{inputenc}
\usepackage{newlfont}
\usepackage{showidx,color}
\usepackage[dvips]{graphicx}
\usepackage{hyperref}
\usepackage[french]{babel}

\begin{document}
\author{M.Pevzner\footnote{This research has been supported by a Marie Curie Fellowship of the
European Community programme "Improving the Human Research
Potential and the Social-Economic Knowledge Base" under contract
number HPMF-CT-2002-01832.}\;\footnote {Nouvelle adresse: UMR 6056
du CNRS, Université de Reims, Moulin de la Housse, B.P. 1039
F-51687 Reims.}
,$\:$ Ch.Torossian\footnote{UMR 8553 du CNRS.} \\
\\ {\footnotesize D\'epartement de Math\'ematiques et Applications- UMR 8553 du CNRS}\\
{\footnotesize \'Ecole Normale Sup\'erieure de Paris}\\
{\footnotesize 45 rue d'Ulm, F-75230 Paris, Cedex 05}\\
{\footnotesize E-mail: \texttt{pevzner@ens.fr,
charles.torossian@ens.fr}}}
\title{Isomorphisme de Duflo et la cohomologie tangentielle }
\date{}
\maketitle
 \pagestyle{myheadings}
\markboth{M.Pevzner, Ch.Torossian}{Isomorphisme de Duflo et la
cohomologie tangentielle }
% \addtolength{\topmargin}{-1.5cm}
\renewcommand{\theequation}{\mbox{\arabic{section}.\arabic{equation}}}
\newtheorem{thm}{Th\'eor\`eme}[section]
\newtheorem{lem}[thm]{Lemme}
\newtheorem{englth}[thm]{Theorem}
\newtheorem{pro}[thm]{Proposition}
\newtheorem{cor}[thm]{Corollary}
\newtheorem{de}[thm]{Definition}
\newtheorem{defn}[thm]{Definition}

          \newtheorem{dfn}{Definition}[section]
          \newtheorem{rmk}{Remark}[section]
          \newtheorem{prop}{Proposition}[section]
          \newtheorem{exs}{Examples}[section]
          \newcommand{\D}{\mbox{$\cal D$}}
          \newcommand{\BR}{\mbox{$\Bbb R$}}
          \newcommand{\BC}{\mbox{$\Bbb C$}}
          \newcommand{\Levi}{\mbox{$\cal L$}}
          \newcommand{\Arg}{\mbox{\rm Arg}}
          \newcommand{\En}{\mbox{$E[[\nu]]$}}
          \newcommand{\Lg}{\mbox{$\frak g$}}
          \newcommand{\Lk}{\mbox{$\frak k$}}
          \newcommand{\La}{\mbox{$\frak a$}}
          \newcommand{\Lh}{\mbox{$\frak h$}}
          \newcommand{\Lm}{\mbox{$\frak m$}}
          \newcommand{\Ln}{\mbox{$\frak n$}}
          \newcommand{\Lt}{\mbox{$\frak t$}}
          \newcommand{\Ll}{\mbox{$\frak l$}}
          \newcommand{\Lz}{\mbox{$\frak z$}}
          \newcommand{\LH}{\mbox{$\frak H$}}
          \newcommand{\Lu}{\mbox{$\frak u$}}
          \newcommand{\Lb}{\mbox{$\frak b$}}
          \newcommand{\Lc}{\mbox{$\frak c$}}
          \newcommand{\Ls}{\mbox{$\frak s$}}
          \newcommand{\Lr}{\mbox{$\frak r$}}
          \newcommand{\Lq}{\mbox{$\frak q$}}
          \newcommand{\Le}{\mbox{$\frak e$}}
          \newcommand{\Lhr}{\mbox{${\frak h}_r$}}
          \newcommand{\Lqr}{\mbox{${\frak p}_r$}}
          \newcommand{\Lhs}{\mbox{${\frak h}_s$}}
          \newcommand{\Lqs}{\mbox{${\frak p}_s$}}
          \newcommand{\Lsu}{\mbox{$\frak{su}$}}
          \newcommand{\Lso}{\mbox{$\frak{so}$}}
          \newcommand{\Lsp}{\mbox{$\frak{sp}$}}
          \newcommand{\Lspin}{\mbox{$\frak{spin}$}}
          \newcommand{\Lgl}{\mbox{$\frak{gl}$}}
          \newcommand{\ad}{\mbox{ad}}
          \newcommand{\adk}{\mbox{$\mbox{ad}_{\Lk}$}}
          \newcommand{\adh}{\mbox{$\mbox{ad}_{\Lh}$}}
          \newcommand{\adx}{\mbox{$\mbox{ad}_{\xi}$}}
          \newcommand{\pr}{\mbox{pr}}
          \newcommand{\Pf}{{\em Proof}. }
          \newcommand{\EPf}{\hfill$\Box$}
          \newcommand{\holm}{\frak{hol}}
          \newcommand{\rad}{\mbox{rad}}
          \newcommand{\bxi}{\mbox{$\bar{X}_i$}}
          \newcommand{\bxj}{\mbox{$\bar{X}_j$}}
          \newcommand{\bxk}{\mbox{$\bar{X}_k$}}
          \newcommand{\byi}{\mbox{$\bar{Y}_i$}}
          \newcommand{\byj}{\mbox{$\bar{Y}_j$}}
          \newcommand{\byk}{\mbox{$\bar{Y}_k$}}
          \newcommand{\tOmega}{\mbox{$\tilde{\Omega}$}}
          \newcommand{\tnabla}{\mbox{$\tilde{\nabla}$}}
          \newcommand{\hOmega}{\mbox{$\hat{\Omega}$}}
          \newcommand{\OO}{{\mbox{${\cal O}$}}}

%          \renewcommand{\theenumi}{\(\alph{enumi}\)}
%\setcounter{section}{-1}
%%% d/but
\let\Dbl\Bbb
\let\BBox\Box
\def\Box{$\BBox$}
%\font\goth=eufm10
\def\g{\mathfrak{g}}
\def\q{\mathfrak{q}}
\def\k{\mathfrak{k}}
\def\t{\mathfrak{t}}
\def\s{\mathfrak{s}}
\def\u{\mathfrak{u}}
\def\z{\mathfrak{z}}
\def\a{\mathfrak{a}}
\def\n{\mathfrak{n}}
\def\p{\mathfrak{p}}
\def\P{\mathfrak{P}}
\def\h{\mathfrak{h}}
\def\gl{\mathfrak{gl}}
\def\sl{\mathfrak{sl}}
\def\e{\mathfrak{e}} \def\l{\mathfrak{l}}
\def\V{\mathfrak{V}} \def\W{\mathfrak{W}}

\def\uple#1#2{#1_1,\ldots ,{#1}_{#2}}
\def\C{\mathbb{C}}
\def\R{\mathbb{R}}
\def\H{\mathbb{H}}
\def\N{\mathbb{N}}
\def\Z{\mathbb{Z}}
\def\qed{\hfill\quad\raise -2pt \hbox{\vrule\vbox to 10pt
{\hrule width 8pt
\vfill\hrule}\vrule}\newline}
\def\gl{\lambda}
\def\C{{\mathbb C}}
\def\Hil{{\mathcal H}}
\def\Cci{C_c^\infty}
\def\cF{{\mathcal F}}
\def\Fougd{\cF_\gd}
\def\Hom{{\rm Hom}}
\def\fa{{\frak a}}
\def\iC{{\scriptstyle\C}}
\def\fadc{\fa_{\iC}^*}
\def\gd{\delta}
\def\fadp{\fa^{*+}}
\def\HIL{{\frak H}}
\def\Fou{{\cF}}
\def\gf{\varphi}
\def\after{{\scriptstyle\circ}}
\def\inp#1#2{\langle #1,#2 \rangle}
\def\ge{\epsilon}

\def\Det{\mathop{\rm Det}\nolimits}
\def\Ad{\mathop{\rm Ad}\nolimits}
\def\id{\mathop{\rm id}\nolimits}
\def\Re{\mathop{\rm Re}\nolimits}
\def\Exp{\mathop{\rm Exp}\nolimits}

\def\tr{\mathop{\rm tr}\nolimits}
\def\Lie{\mathop{\rm Lie}\nolimits}

\def\spur{\mathop{\rm spur}\nolimits}

\def\div{\mathop{\rm div}\nolimits}

\ifx\optionkeymacros\undefined\else\endinput\fi

\def\mod{\mathop{\rm mod}\nolimits}

\def\sgn{\mathop{\rm{sgn}}\nolimits}

\def\Tr{\mathop{\rm Tr}\nolimits}

\def\rank{\mathop{\rm rank}\nolimits}

\def\diag{\mathop{\rm diag}\nolimits}

\def\Int{\mathop{\rm int}\nolimits}
\def\carre{\mathbin{\hbox{\vrule\vbox to 4pt
    {\hrule width 4pt \vfill\hrule}\vrule}}}

\def\Aut{\mathop{\rm Aut}\nolimits}

\def\ch{\mathop{\rm ch}\nolimits}

\def\tanhyp{\mathop{\rm th}\nolimits}
\def\adnu{\mathop{\rm ad}_{\star_{\nu}}\nolimits}
\def\Exp{\mathop{\rm Exp}\nolimits}

\def\ad{\mathop{\rm ad}\nolimits}

\def\Im{\mathop{\rm Im}\nolimits}

\def\Ad{\mathop{\rm Ad}\nolimits}

\def\dim{\mathop{\rm dim}\nolimits}

\def\Sym{\mathop{\rm Sym}\nolimits}
\def\rm{\mathrm}
\abstract{In the present note we show that the Duflo isomorphism
extends to an isomorphism of associative algebras of tangential
cohomologies. This result confirms the B.Shoikhet's conjecture
\cite{[Sh1],[Sh2]}.} \vskip 10pt \noindent{\bf 2000 Mathematics
Subject classification} : 16S80, 53D17, 53D55, 22E30, 22E60.
\tableofcontents
\section*{Introduction}
Soit $X$ une vari\'et\'e de Poisson quelconque. Le th\'eor\`eme de
formalit\'e de M. Kon\-tsevich \cite{[K97]} Th.6.4 montre
l'existence d'un $L_{\infty}$-quasi-isomor\-phisme entre deux
alg\`ebres de Lie diff\'erentielles gradu\'ees naturellement
associ\'ees \`a la vari\'et\'e $X$. Plus pr\'ecisement, ces deux
alg\`ebres $\g_1$ et $\g_2$ sont respectivement, celle des
polychamps de vecteurs sur $X$ munie de la diff\'erentielle nulle
et du crochet de Schouten-Nijenhuis et celle des op\'erateurs
polydiff\'erentiels munie de la diff\'eren\-tielle de Hochschild
et du crochet de Gerstenhaber.\\

Lorsque la vari\'et\'e $X$ est le dual d'une alg\`ebre de Lie
r\'eelle de dimension finie $\g$ il existe sur $X$  une structure
de Poisson canonique induite par le crochet de Kirillov-Kostant.
Dans ce cas particulier la structure de Poisson est lin\'eaire
(les coefficients du bivecteur correspondant sont des fonctions
lin\'eaires). Le th\'eor\`eme de la formalit\'e implique alors
l'existence d'un isomorphisme d'alg\`ebres entre l'ensemble
$S(\g)^{\g}$ des polyn\^omes $\g$-invariants sur $\g^*$ et le
centre $Z(\g)$ de l'alg\`ebre enveloppante de $\g$.

Cet isomorphisme est donn\'e par la diff\'erentielle du
$L_{\infty}$-quasi-isomorphisme de Kontsevich restreinte aux
0-cohomologies tangentes des alg\`ebres $\g_1$ et $\g_2$
associ\'ees \`a $X=\g^*$.

De plus cet isomorphisme co\"\i ncide avec l'isomorphisme de
\mbox{Duflo \cite{[K97]}.}\\

Dans la pr\'esente note nous montrons que cet isomorphisme se
prolonge en cohomologie tangentielle de plus haut degr\'e. Ce
r\'esultat confirme une conjecture de B.Shoikhet
\cite{[Sh1],[Sh2]}. \setcounter{equation}{0}
\section{Rappels et Notations}
\subsection{Alg\`ebre des polychamps de vecteurs}
 Soit $X$ une vari\'et\'e  de classe $C^{\infty}$. On
lui associe
 deux alg\`ebres de Lie diff\'erentielles gradu\'ees. La
 premi\`ere alg\`ebre de Lie diff\'erentielle gradu\'ee~
 $\g_1=T_{\rm{poly}}(X)$ est l'alg\`ebre gradu\'ee des
 polychamps de vecteurs sur $X$:
 $$
 T^n_{\rm{poly}}(X):=\Gamma(X,\Lambda^{n+1}TX),\quad n\geq-1
 $$
munie du crochet de Schouten-Nijenhuis $[\;,\:]_{SN}$ et de la
diff\'erentielle $d:=0$.

Rappelons tout d'abord que le crochet de Schouten-Nijenhuis est
donn\'e pour tous $k,l\geq 0$ , $\xi_i, \eta_j \in \Gamma(X,TX)$
par (\cite{[Kosz]}):
\begin{eqnarray*}
[\xi_0\wedge\dots\wedge\xi_k,\eta_0\wedge\dots\wedge\eta_l]_{SN}=
\sum_{i=0}^k\sum_{j=0}^l(-1)^{i+j}[\xi_i,\eta_j]
\wedge&&\\\xi_0\wedge\dots\wedge\xi_{i-1}\wedge\xi_{i+1}\wedge\ldots\wedge\xi_k
\wedge\eta_0\wedge\dots\wedge\eta_{j-1}\wedge\eta_{j+1}\wedge\dots\wedge\eta_l&&
\end{eqnarray*}
Et pour $k\geq 0$ et $h\in\Gamma(X,{\cal
O}_X),\:\xi_i\in\Gamma(X,TX)$:
$$
[\xi_0\wedge\dots\wedge\xi_k,h]_{SN}=
\sum_{i=0}^k(-1)^{i}\xi_i(h)\cdot\left(\xi_0\wedge\dots\wedge\xi_{i-1}\wedge\xi_{i+1}\wedge\ldots\wedge\xi_k\right).
$$
Le symbole $[\xi_i,\eta_j]$ d\'esigne le crochet standard des
champs des vecteurs, {\it i.e.} la d\'eriv\'ee de Lie
$L_{\xi_i}(\eta_j)$.\\

Pr\'ecisons ici que nous adoptons la convention de \cite{[MT]}
l'identification des poly-champs de vecteurs avec les tenseurs
anti-sym\'etriques. Les produits $\xi_1\wedge\ldots\wedge \xi_k$
s'identifient avec
$${1\over k!}\sum_{\sigma\in S_k}\sgn(\sigma)\xi_{\sigma_1}\otimes\cdots\otimes
\xi_{\sigma_k}.$$Cette convention est diff\'erente de celle de
\cite{[K97]} \S 6.3 mais sera compens\'ee par une autre convention dans la d\'efinition des poids.\\

Nous renvoyons le lecteur vers \cite{[K97]} \S4.6.1 et
\cite{[AMM]} IV.2 pour une autre interpr\'etation du crochet de
Schouten-Nijenhuis en termes des varia\-bles impaires.
\subsection{Alg\`ebre  des op\'erateurs polydiff\'erentiels}
 La deuxi\`eme alg\`ebre de Lie diff\'erentielle
gradu\'ee associ\'ee \`a $X$ est celle des op\'erateurs
polydiff\'erentiels  $\g_2=D_{\rm{poly}}(X)$ vue comme une
sous-alg\`ebre du complexe de Hochschild d\'ecal\'e de l'alg\`ebre
des fonctions sur $X$.

On d\'efinit sur $D_{\rm{poly}}(X)$ une graduation donn\'ee par
$\vert A\vert=m-1$ o\`u $A\in D_{\rm{poly}}(X)$ est un op\'erateur
$m-$diff\'erentiel.\\

La composition de deux op\'erateurs $A_1\in
D^{m_1}_{\rm{poly}}(X)$ et $A_2\in D^{m_2}_{\rm{poly}}(X)$
s'\'ecrit pour $f_i\in {\cal O}_X$:
\begin{eqnarray}
(A_1\circ
A_2)(f_1,\ldots,f_{m_1+m_2-1})&=&\sum_{j=1}^{m_1}(-1)^{(m_2-1)(j-1)}
A_1(f_1,\ldots, f_{j-1},\nonumber\\
&& \hspace{-2cm}A_2(f_j,\ldots,f_{j+m_2-1}),f_{j+m_2},
\ldots,f_{m_1+m_2-1}).\nonumber
\end{eqnarray}
 Cette
op\'eration de composition permet de d\'efinir le crochet de
Gerstenhaber comme suit:
$$
[A_1,A_2]_G:=A_1\circ A_2-(-1)^{\vert A_1\vert\vert
A_2\vert}A_2\circ A_1.
$$
La diff\'erentielle dans $D_{\rm{poly}}(X)$ s'\'ecrit alors
$$
dA=-[\mu,A]_G,
$$
o\`u $\mu$ est l'op\'erateur bi-diff\'erentiel de multiplication
des fonctions:
$$\mu(f_1,f_2)=f_1f_2.$$

Notons que cette diff\'erentielle est li\'ee \`a celle de
Hochschild $d_H$ par la relation $$d_H(A)=(-1)^{\vert
A\vert+1}dA.$$

Ce choix de signe fait de $(D_{\rm{poly}}(X),d,[\:,\:]_G)$ une
alg\`ebre de Lie diff\'erentielle gradu\'ee.

\subsection{Th\'eor\`eme de formalit\'e de M. Kontsevich}
L'application ${\cal U}_1^{(0)}:T_{\rm{poly}}\mapsto
D_{\rm{poly}}$ donn\'ee par
\begin{equation}\label{uzero}
{\cal U}_1^{(0)}:(\xi_0,\dots,\xi_n)\mapsto
\left(f_0\otimes\dots\otimes f_n\to\sum_{\sigma\in
S_{n+1}}\frac{\sgn(\sigma)}{(n+1)!}\prod_{i=0}^n\xi_{\sigma(i)}(f_i)\right)
\end{equation}
pour $n\geq 0$ et par $f\mapsto (1\to f)$ pour $f\in\Gamma(X,{\cal
O}_X)$ est un quasi-isomorphisme des complexes. C'est une version
du th\'eor\`eme de Kostant-Hochschild-Rosenberg (\cite{[K97]}
\S4.6.1.1).

\vskip 15pt
 Rappelons bri\`evement la construction du
$L_{\infty}$-quasi-isomorphisme de Kontsevich en suivant la
pr\'esentation de \cite{[K97],[AMM],[MT]}.\\

Consid\'erons les alg\`ebres d\'ecal\'ees $\g_1[1]$ et $\g_2[1]$
comme des vari\'et\'es formelles gradu\'ees, point\'ees ({\it cf}
\cite{[K97]}, \cite{[AMM]} III.2). Chacune des cog\`ebres sans
co-unit\'e
$$
S^+\big(\g_i[1]\big)=\bigoplus_{n\geq 0}S^n\big(\g_i[1]\big),
\quad i=1,2
$$
 poss\`ede une cod\'eri\-vation $Q^i$ de degr\'e 1 d\'efinie
par la structure d'alg\`ebre de Lie diff\'erentielle
gradu\'ee\footnote{D'apr\`es \cite{[AMM]} on doit remplacer le
crochet de Schouten par l'oppos\'e du crochet pris dans l'ordre
inverse. Ce crochet co\"\i ncide avec le crochet de Schouten
modulo un signe moins lorsque deux \'el\'ements impairs sont en
jeu. Remarquons qu'alors le c\oe fficient de Taylor $Q_2^1$
d\'efini sur $S(\g_1[1])$ vaut alors pour
$\gamma_1=\xi_1\wedge\ldots \wedge \xi_{k_1}$ et
$\gamma_2=\eta_1\wedge\ldots \wedge \eta_{k_2}$ :
$$Q_2^1(\gamma_1.\gamma_2)=\sum_{\substack{1\leq r \leq k_1\atop
1\leq s\leq k_2}}(-1)^{r+s+k_1-1}[\xi_r,\eta_s]\wedge
 \xi_1\wedge \ldots \widehat{\xi_r}\wedge \ldots \xi_{k_1}
\wedge\eta_1\wedge\ldots \widehat{\eta_s}\wedge\ldots
\eta_{k_2},$$ ce qui est conforme \`a la formule de \cite{[K97]}
\S4.6.1.} de $\g_i$ et telle que $[Q^i,Q^i]=0$. Le th\'eor\`eme de
la formalit\'e peut \^etre \'enonc\'e alors de la fa\c con
suivante.
\begin{thm}\label{formal}
Il existe un $L_{\infty}$-quasi-isomorphisme entre les
vari\'et\'es formelles gradu\'ees point\'ees $\g_1[1]$ et
$\g_2[1]$, i.e. un morphisme de cog\`ebres
$$
{\cal U}:S^+\big(\g_1[1]\big)\to S^+\big(\g_2[1]\big)
$$
tel que
$$
{\cal U}\circ Q^1=Q^2\circ{\cal U}
$$
et tel que la restriction de ${\cal U}$ \`a $\g_1[1]\simeq
S^1(\g_1[1])$ est le quasi-isomorphisme de complexe de cocha\^\i
nes ${\cal U}_1^{(0)}$ donn\'e par (\ref{uzero}).
\end{thm}

 Cette construction est bas\'ee sur la
 r\'ealisation explicite du $L_{\infty}$-quasi-isomorphisme $\cal U$ dans le cas plat, dont nous allons
 \`a pr\'esent rappeler la construction.
Nous suivons la pr\'esentation de  \cite{[K97]} chapitres 5 et 6.\\

Pla\c cons nous dans le cas o\`u la vari\'et\'e $X$ est l'espace
vectoriel $\Bbb R^d$. D'apr\`es la propri\'et\'e universelle des
cog\`ebres,
 cocommutatives le
 $L_{\infty}$-quasi-isomorphisme $\cal U$ est enti\`erement
 d\'etermin\'e par ses "c\oe fficients de Taylor"
 $$
 {\cal U}_k:S^k(\g_1[1])\to\g_2[1],
$$
 avec $k\geq 1$ obtenus en composant $\cal U$ avec
 la projection canonique $\pi:S^+(\g_2)\to\g_2$. On notera
 $\overline{\mathcal{U}}$ cette composition.\\

Ces c\oe fficients de Taylor ${\cal U}_n$  sont  d\'efinis \`a
l'aide des graphes de Kontsevich et de leurs poids dont nous
rappelons bri\`evement la construction afin de rendre le texte
plus autonome et de fixer quelques conventions.

 \vskip 15pt
 Soit $G_{n,m}$ l'ensemble de graphes \'etiquet\'es,
orient\'es ayant $n$ sommets du premier type (sommets a\'eriens)
et $m$ sommets du deuxi\`eme type (sommets terrestres) tels que:\\

1. Toutes les ar\^etes partent des sommets du premier type.

2. Le but d'une ar\^ete est diff\'erent  de sa source (pas de
boucles).

3. Il n'y a pas d'ar\^etes multiples (m\^eme source, m\^eme
but).\\

\noindent On dira qu'un tel graphe est admissible. Par
\'etiquetage d'un graphe admissible $\Gamma$ on entend un ordre
total sur l'ensemble $E_{\Gamma}$ des ar\^etes de $\Gamma$
compatible avec l'ordre de l'ensemble des sommets.\\

Soit  $\Gamma\in G_{n,m}$ un graphe admissible \'etiquet\'e, on
note  $s_k$  le nombre d'ar\^etes partant du sommet du premier
type ayant le num\'ero $k$. A tout $n$-uplet
$(\alpha_1,\dots,\alpha_n)$ de polychamps de vecteurs sur $X$ tels
que pour tout $k=1,\dots,n$ l'\'el\'ement $\alpha_k$ soit un
$s_k$-champ de vecteurs, on associe, suivant \cite{[K97]} \S 6.3
un op\'erateur $m$-diff\'erentiel
$$
B_{\Gamma}(\alpha_1\otimes\dots\otimes\alpha_n),$$  construit de
la fa\c con suivante~: on d\'esigne par
$\{e_k^1,\ldots,e_k^{s_k}\}$ le sous-ensemble ordonn\'e de
$E_\Gamma$ des ar\^etes partant du sommet a\'erien $k$. A toute
application $I:E_\Gamma\to \{1,\ldots, d\}$ et \`a tout sommet $x$
du graphe $\Gamma$ (de type a\'erien ou terrestre) on associe
l'op\'erateur diff\'erentiel \`a c\oe fficient constant~:
$$D_{I(x)}=\prod_{e=(-,x)}\partial_{I(e)},$$
o\`u pour tout $i\in\{1,\ldots,d\}$ on d\'esigne par $\partial_i$
l'op\'erateur de d\'erivation partielle par rapport \`a la
$i$-\`eme variable. Le produit est pris pour toutes les ar\^etes
qui arrivent au sommet $x$.\\

\noindent On d\'esigne par $\alpha_k^I$  le c\oe fficient (suivant
la convention sur le produit ext\'erieur)~:
\begin{eqnarray*}
\alpha_k^I=\alpha_k^{I(e_k^1)\cdots I(e_k^{s_k})} &=&
\langle\alpha_k,\,dx_{I(e_k^1)}\wedge\dots\wedge dx_{I(e_k^{s_k})}\rangle\\
        &=&\langle\alpha_k,dx_{I(e_k^1)}\otimes\cdots\otimes dx_{I(e_k^{s_k})}\rangle.
\end{eqnarray*}
         On pose alors~:
$$ B_\Gamma(\alpha_1\otimes\cdots\otimes \alpha_n)(f_1\otimes\cdots\otimes f_m)=\sum_{I:E_\Gamma\to \{1,\ldots, d\}}\prod_{k=1}^n D_{I(k)}\alpha_k^I
        \prod_{l=1}^mD_{I(\overline l)}f_l.
$$
\smallskip
Le c\oe fficient de Taylor ${\cal U}_n$ est donn\'e par la
formule~:
\begin{equation}\label{Un}
{\cal U}_n(\alpha_1,\ldots ,\alpha_n)=\sum_{\Gamma\in G_{n,m}}
        w_\Gamma B_\Gamma(\alpha_1\otimes\cdots\otimes\alpha_n),
\end{equation}
o\`u la somme porte sur les graphes admissibles $\Gamma$ pour
lesquels l'op\'erateur
$B_\Gamma(\alpha_1\otimes\cdots\otimes\alpha_n)$ est bien d\'efini
et l'entier $m$ est reli\'e \`a $n$ et aux $\alpha_j$ par la
formule
\begin{equation}\label{M} m-2=\sum_{k=1}^ns_k-2n.\end{equation}\\

Le c\oe fficient $w_\Gamma$ est un certain poids associ\'e \`a
chaque graphe $\Gamma$. Alors ${\cal
U}_n(\alpha_1,\dots,\alpha_n)$ est un op\'erateur
$m-$diff\'erentiel.\\

 Le poids $w_\Gamma$ est nul sauf si le
nombre d'ar\^etes $|E_\Gamma|$ du graphe $\Gamma$ est
pr\'ecis\'ement \'egal \`a $2n+m-2$. Il s'obtient en int\'egrant
une forme ferm\'ee $\Omega_\Gamma$ de degr\'e $|E_\Gamma|$ sur une
composante connexe de la compactification de Fulton-McPherson d'un
espace de configurations qui est pr\'ecis\'ement de dimension
$2n+m-2$ (voir \cite {[FM]} et \cite {[K97]} \S 5). Ce poids
d\'epend lui aussi d'un ordre sur l'ensemble des ar\^etes, mais le
produit $w_\Gamma\cdot B_\Gamma$ n'en d\'epend plus.\\

 Plus pr\'ecis\'ement on
d\'esigne par ${\rm{Conf}}_{n,m}$ l'ensemble des $$(\uple pn,\uple
qm)$$o\`u les $p_j$ sont des points distincts appartenant au
demi-plan de Poincar\'e~:
$$\Bbb H_+=\{z\in\Bbb C, \Im z>0\},$$
et o\`u les $q_j$ sont des points distincts sur $\Bbb R$ vu comme
le bord de $\Bbb H_+$. Le groupe~:
$$G=\{z\mapsto az+b \hbox{ avec }(a,b)\in\Bbb R \hbox{ et }a>0\}$$
agit librement sur ${\rm{Conf}}_{n,m}$. Le quotient~:
\begin{equation}\label{cnm}
C_{n,m}={\rm{Conf}}_{n,m}/G
\end{equation}
est une vari\'et\'e de dimension $2n+m-2$. Dans \cite {[K97]} \S
5, M. Kontsevich construit des compactifications
$\overline{C_{n,m}}$ de ces vari\'et\'es de configurations. Ce
sont des vari\'et\'es \`a
coins de dimension $2n+m-2$.\\

\begin{figure}[!h]
\begin{center}

\includegraphics[width=2in]{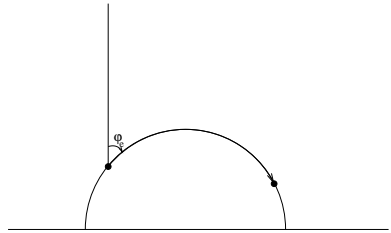}
  \caption{Fonction d'angle}\label{angle}
  \end{center}
\end{figure}

Pour tout graphe $\Gamma\in G_{n,m}$ on d\'efinit une fonction
d'angle~:
\begin{equation}\label{weight}
\Phi_\Gamma:\overline C_{n,m}\longrightarrow (\Bbb R/2\pi\Bbb
Z)^{|E_\Gamma|}
\end{equation}
de la fa\c con suivante~: on trace le graphe dans $\overline{\Bbb
H_+}$ en reliant les sommets par des g\'eod\'esiques pour la
m\'etrique hyperbolique, et \`a chaque ar\^ete $e=(p,q)$ on
associe l'angle $\varphi_e={\mathrm Arg}\left(\frac{q-p}{q-\bar
p}\right)$ que fait la demi-droite verticale issue de $p$ avec
l'ar\^ete $e$ (voir figure [\ref{angle}]).

En choisissant un ordre sur les ar\^etes ceci d\'efinit
$\Phi_\Gamma$ sur $C_{n,m}$ et on v\'erifie que cette application
se prolonge \`a la compactification. Soit $\Omega_\Gamma$ la forme
diff\'erentielle $\Phi_\Gamma^*(dv)$ sur $\overline{C}_{n,m}$ o\`u
$dv$ est la forme volume normalis\'ee sur $(\Bbb R/2\pi\Bbb
Z)^{|E_\Gamma|}$. Soit $\overline{C}^+_{n,m}$ la composante
connexe de $\overline{C}_{n,m}$ o\`u les $\uple qm$ sont rang\'es
par ordre croissant. Les orientations naturelles du demi-plan
$\Bbb H_+$ et de $\Bbb R$ d\'efinissent une orientation de
${\rm{Conf}}_{n,m}^+$, et par passage au quotient une orientation
naturelle de $\overline{C}_{n,m}^+$, car l'action du groupe $G$
pr\'eserve l'orientation. On d\'efinit alors le poids $w_\Gamma$
par~:
\begin{equation}
w_\Gamma=\int_{\overline C_{n,m}^+}\Omega_\Gamma.
\end{equation}
{\it Remarque}~: Ce poids est un peu diff\'erent du poids d\'efini
par M. Kontsevich dans \cite {[K97]} \S\ 6.2~: nous ne multiplions
pas l'int\'egrale par le facteur $\bigl(\prod_{k=1}^n{1\over
s_k!}\bigr)$. Cette convention est compens\'ee par celle sur les sur le produit ext\'erieur.\\
%
%\subsection{ Permutation des ar\^etes}
%

 Soit $\Gamma$ un graphe admissible dans $G_{n,m}$. Le groupe
$S_{s_1}\times\cdots\times S_{s_n}$, produit des groupes de
permutations des ar\^etes attach\'es \`a chaque sommet, agit
naturellement sur $\Gamma$ par permutation de l'\'etiquetage des
ar\^etes. Il est clair que l'on a~:
\begin{eqnarray*}
B_{\sigma.\Gamma}&=&\varepsilon(\sigma) B_\Gamma \\
w_{\sigma.\Gamma}&=&\varepsilon(\sigma)w_\Gamma,
\end{eqnarray*}
de sorte que le produit $w_\Gamma\cdot B_\Gamma$ ne d\'epend pas
de l'\'etiquetage.

 \vskip 10pt
% La
%diff\'erentielle ${\cal U}_{\hbar\gamma}'$ du
%$L_{\infty}$-quasi-isomorphisme $\cal U$ au $\Lm$-point
%$\hbar\gamma$ peut \^etre d\'ecrite explicitement de la fa\c con
%suivante. On identifie l'espace tangent $T_{\hbar\gamma}(\g_1[1])$
%\`a $\g_1[1]\hat\otimes\Lm$. Par lin\'earisation un champ de
%vecteurs impair $Q^1$ s'annulant en $\hbar\gamma$ d\'efinit un
%champ de vecteurs impair $Q^{\hbar\gamma}$ sur
%$T_{\hbar\gamma}(\g_1[1])$ par
%\begin{equation}
%Q^{\hbar\gamma}(\hbar\delta)=\sum_{n>0}\frac{\hbar^{n+1}}{n!}Q_{n+1}^1(\delta.\gamma^{.n}).
%\end{equation}
%On obtient en occurrence
%\begin{equation}
%Q^{\hbar\gamma}(\hbar\delta)=Q_2^1(\hbar\gamma.\hbar\gamma)=\hbar^2[\delta,\gamma]_{SN}.
%\end{equation}
%De m\^eme pour le champ de vecteurs impairs $Q^2$ on obtient un
%champ de vecteurs $Q^{\hbar\tilde\gamma}$ de carr\'e nul sur
%$T_{\hbar\tilde\gamma}(\g_2[1])$ qui s'\'ecrit ({\it cf}
%\cite{[AMM]} II.4 et IV.1):
%\begin{equation}
%Q^{\hbar\tilde\gamma}(\hbar\delta)=\hbar[\delta,\star]_{G}.
%\end{equation}
%Ainsi les deux espaces tangents $T_{\hbar\gamma}(\g_1[1])$ et
%$T_{\hbar\tilde\gamma}(\g_2[1])$ sont des complexes de cocha\^\i
%nes et la d\'eriv\'ee ${\cal U}_{\hbar\gamma}'$ du
%$L_{\infty}$-quasi-isomorphisme $\cal U$ est un
%quasi-isomorphismes de ces complexes. Elle s'exprime alors par la
%formule:
%\begin{equation}\label{uderiv}
%{\cal U}_{\hbar\gamma}'(\delta)=\sum_{n\geq
%0}\frac{\hbar^n}{n!}{\cal U}_{n+1}(\delta.\gamma^{.n})
%\end{equation}
%pour tout $\delta\in\g_1[1]$.
\subsection{Quasi-isomorphisme tangentiel et la formule
d'homotopie}

Soit $\gamma\in T_{\rm{poly}}(\Bbb R^d)[1]$ un $2-$champ de
vecteurs tel que $\gamma[-1]$ v\'erifie l'\'equation de
Maurer-Cartan dans $T_{\rm{poly}}(\Bbb R^d)$. C'est donc un
$2-$champ de vecteurs de Poisson. Le
$L_{\infty}$-quasi-isomorphisme introduit par M. Kontsevich permet
de construire, \`a partir de $\hbar \gamma$, un
 star-produit $\star_{\hbar\gamma}$ de la mani\`ere suivante:
 \begin{equation}\label{starproduit}
 \star_{\hbar\gamma}=\mu +\overline{\mathcal{U}}(\hbar \gamma)=\mu +
 \sum_{n\geq 1} \frac{\hbar^n}{n!}\mathcal{U}_n(\gamma, \ldots
 ,\gamma),
 \end{equation}
o\`u $\hbar$ est un param\`etre formel. Nous allons \'etudier en
d\'etail les propri\'et\'es de la d\'eriv\'ee du
$L_{\infty}$-quasi-isomor\-phisme $\cal U$ au point $\hbar\gamma$.
\\

Pour cela nous suivons la pr\'esentation de \cite{[Mosh]}.\\

Nous avons rappel\'e  en (\ref{cnm}) la d\'efinition des
compactifications des espaces de confi\-gurations
$\overline{C_{n,m}}$. Comme une variation de ces objets Mochizuki
consid\`ere les espaces $X_{n,m}^{\ell}$ d\'efinis pour tout
entier positif $\ell$ de fa\c
con suivante.\\

Soit $\phi$ l'isomorphisme naturel entre $C_{\ell+1,0}^{\Bbb R}$
espace de configurations r\'eelles et l'int\'erieur du
$(\ell-1)$-simplexe standard
 ${\overset{\circ}{\Delta}}\, ^{\ell-1}$, donn\'e par $\phi(x)=(p_0(x)=0,\dots,p_{\ell}(x)=1)$.
On obtient ainsi l'application $\psi:{\overset{\circ}{\Delta}}\,
^{\ell}\simeq{\overset{\circ}{\Delta}}\, ^{\ell-1}\times\Bbb
R_+\to
C_{\ell+1,0}$ d\'efinie par $\psi(x,t)=(p_0(x)+it, \ldots, p_{\ell}(x)+it)$.\\

Notons alors
$\overset{\circ}{X}\,_{n,m}^{\ell}=C_{n,m}^+\times_{C_{\ell+1,0}}{\overset{\circ}{\Delta}}\,
^{\ell}$ et $X_{n,m}^{\ell}$ son adh\'erence dans
$\overline{C_{n,m}}$. On obtient ainsi une vari\'et\'e \`a coin de
dimension
$2n+m-2-\ell$.\\

On g\'en\'eralise  la notion d'un graphe admissible. On dira qu'un
graphe admissible $\Gamma$ est de type $(n,m,e)$ si le cardinal de
l'ensemble de ces points a\'eriens $\sharp V_{\Gamma}^1=n$, celui
de l'ensemble des points terrestres $\sharp
V_{\Gamma}^2=m$ et le nombre d'ar\^etes est \'egal \`a $e$.\\

Un graphe $\Gamma$ est dit $\ell$-admissible ($\ell\in\Bbb N$) si
$\Gamma$ est admissible de type $(n,m,e)$ et
$V_{\Gamma}^1=V_{\Gamma}^{1,1}\sqcup V_{\Gamma}^{1,2}$ tels que\\

a. $ V_{\Gamma}^{1,1}=\{1,\dots,\ell+1\}$

b. Pour tout couple $(i,j)\subset V_{\Gamma}^{1,1}$ il n'existe
pas d'ar\^ete les reliant.\\

\noindent On notera $G_{n,m}^{\ell,k}$ l'ensemble des graphes
$\ell$-admissibles de type $(n,m,2n+m-2-\ell-k)$. Si
$e=2n+m-2-\ell$ on notera $G_{n,m}^{\ell,0}=G_{n,m}^{\ell}$.\\

\noindent On d\'efinit alors les poids des graphes
$\ell$-admissibles par
$$
w_{\Gamma}^{\ell}:=\int_{X_{n,m}^{\ell}}\Omega_{\Gamma},
$$
o\`u $\Omega_{\Gamma}$ est la forme diff\'erentielle d\'efinie par
la fonction d'angle (\ref{weight}).\\

Dans \cite{[Mosh]} Mochizuki introduit des applications
\begin{equation}
{\cal U}_{n,m}^{\ell}:=\sum_{\Gamma\in
G_{n,m}^{\ell}}w_{\Gamma}^\ell\cdot
B_{\Gamma}:\bigotimes^n(T_{\rm{poly}}(\Bbb R^d)[1])\to
D_{\rm{poly}}(\Bbb R^d)[1+\ell]
\end{equation}
et note
\begin{equation}
{\cal U}_{n}^{\ell}:=\sum_m{\cal U}_{n,m}^{\ell}.
\end{equation}
Consid\'erons alors les applications
$$
{\cal
U}_{n,m}^{\ell,\gamma}:\bigotimes^{\ell+1}(T_{\rm{poly}}(\Bbb
R^d)[1])\to D_{\rm{poly}}(\Bbb R^d)[1+\ell]
$$
telles que :
$$
{\cal
U}_{n,m}^{\ell,\gamma}(\alpha_1\otimes\dots\otimes\alpha_{\ell+1})=
{\cal
U}_{n,m}^{\ell}(\alpha_1\otimes\dots\otimes\alpha_{\ell+1}\otimes\gamma\otimes\dots\otimes\gamma).
$$

L'application d\'eriv\'ee $d{\overline{\cal U}}:T_{\rm{poly}}(\Bbb
R^d)[1]\to D_{\rm{poly}}(\Bbb R^d)[1][[\hbar]]$ au point $\hbar
\gamma$ est alors d\'efinie par ({\it cf.} formule (\ref{Un}))
\begin{equation}\label {uderiv}
d{\overline{\cal U}}_{\hbar
\gamma}(\delta):=\sum_{n>0}\frac{\hbar^{n-1}}{(n-1)!}{\cal
U}_{n}(\delta.\gamma^{.n}).
\end{equation}
Dans \cite{[Mosh]} sont d\'efinis  les applications
$$ d{\cal U}^{\ell,\gamma}:\bigotimes^{\ell+1}T_{\rm{poly}}[1]\to
D_{\rm{poly}}[1+\ell][[\hbar]]$$ par
\begin{equation}
d{\cal
U}^{\ell,\gamma}(\alpha_1\otimes\dots\otimes\alpha_{\ell+1})=\sum_{n,m}\frac{\hbar^{n-\ell-1}}{(n-\ell-1)!}{\cal
U}_{n,m}^{\ell,\gamma}(\alpha_1\otimes\dots\otimes\alpha_{\ell+1}).
\end{equation}
Notons que l'on a en particulier  $d\overline{{\cal U}}_{\hbar
\gamma}=d{\cal
U}^{0,\gamma}.$\\

\noindent S'il n'y a pas de confusion possible on notera $d{\cal
U}^{\ell}$ au lieu de $d{\cal U}^{\ell,\gamma}$. Ces applications
permettent de montrer l'existence d'une structure
$A_\infty$ tangente \cite{[Mosh]}.\\

 L'op\'erateur de cobord du complexe
des cocha\^\i nes tangentielles de la première alg\`ebre de Lie
diff\'erentielle gradu\'ee $T_{\rm{poly}}(X)[1]$
 est donn\'e par $Q^{\hbar\gamma}=-[\hbar\gamma,-]_{SN}$ qui est une d\'erivation
gradu\'ee pour le produit ext\'erieur $\wedge$ des poly-champs de
vecteurs (graduation standard).

Ce produit ext\'erieur induit donc un produit associatif et
commutatif que l'on notera $\cup$ sur l'espace de la cohomologie
$H_{\hbar\gamma}$ du premier espace tangent.\\

Sur le deuxi\`eme espace tangent on introduit un produit
associatif gradu\'e donn\'e par la formule suivante:
\begin{eqnarray}
\lefteqn{(A_1\cup A_2)(f_1\otimes\dots\otimes
f_{m_1+m_2})=}\nonumber\\&&A_1(f_1\otimes\dots\otimes
f_{m_1})\star_{\hbar\gamma} A_2(f_{m_1+1}\otimes\dots\otimes
f_{m_2})
\end{eqnarray}
pour tout op\'erateur $m_1$-diff\'erentiel $A_1$ et tout
op\'erateur $m_2$-diff\'erentiel $A_2$. Cette op\'eration est
compatible avec le cobord $[-,\star]_G$ du deuxi\`eme complexe des
cocha\^\i nes tangentielles et elle induit donc un cup-produit sur
l'espace de la cohomologie
$H_{\overline{\mathcal{U}}(\hbar\gamma)}$ du deuxi\`eme espace tangent.\\

Le th\'eor\`eme suivant est d\'emontr\'e avec diff\'erents
degr\'es de pr\'ecision dans \cite{[K97]} \S 8 et \cite{[MT]}
Th\'eor\`eme 1.2.\footnote{Nous profitons de l'occasion pour
corriger dans \cite{[MT]} une erreur de signe dans la proposition
4.1 et le th\'eor\`eme 4.6 due \`a une confusion dans
l'utilisation du lemme 4.2. Dans la proposition 4.1 et le
th\'eor\`eme 4.6,  $\alpha$ est un $k_1$-champ de vecteurs,
$\beta$ est un $k_2$-champ de vecteur et $\gamma$ est le 2-tenseur
de Poisson. Mais dans le lemme 4.2 l'entier $k_2$ se réfère à
$\gamma$, donc est \'egal \`a $2$. Le signe $(-1)^{(k_1-1)k_2}$
vaut donc  $1$. Le signe $(-1)^{k_1(k_2-1)}$ doit \^etre
remplac\'e par $(-1)^{k_1}$, car l'argument \`a la fin de la
d\'emonstration du lemme 4.2 est d\'efaillant et doit \^etre
corrig\'e par le suivant: on \'echange les positions $1$
associ\'ee \`a $\alpha$ et $2$ associ\'ee \`a $\beta$ ce qui fait
appara\^ \i tre un signe $(-1)^{k_1k_2}$ d\^u au poids, puis on
\'echange \`a nouveau les positions apr\`es contraction cette fois
entre $[\beta, \gamma]$ et $\alpha$, ce qui fait appara\^ \i tre
un signe $(-1)^{(k_2+1)k_1}$ et fournit le signe $(-1)^{k_1}$.
Signalons aussi que dans l'expression du crochet de Schouten
modif\'e apr\`es le lemme 4.2, le signe $(-1)^{k_1k_2}$ doit
\^etre remplac\'e par $(-1)^{(k_1-1)(k_2-1)}$ ce qui entraîne  une
modification de signe dans la formule $Q_2^1(\gamma_1.\gamma_2)$
quelques lignes plus loin, on doit remplacer $(-1)^{k_2}$ par
$(-1)^{k_1-1}$. Toutes ces modifications sont mineures et sans
cons\'equences sur le reste de l'article.}

Dans \cite{[Mosh]} on trouvera une preuve de l'existence d'une
$A_\infty$ structure tangente.
\begin{thm}\label{homotopie}
Soit $X=\Bbb R^d$ et $\cal U$ le $L_{\infty}$-quasi-isomorphisme
donn\'e explicitement dans \cite{[K97]} \S 6.4.

La d\'eriv\'ee $d{\overline{\cal U}}_{\hbar\gamma}$ induit un
isomorphisme d'alg\`ebres de l'espace de la cohomologie
$H_{\hbar\gamma}$ de l'espace tangent $T_{\hbar\gamma}(\g_1[1])$
sur l'espace de la cohomologie
$H_{{\overline{\mathcal{U}}(\hbar\gamma)}}$ de l'espace tangent
$T_{{\overline{\mathcal{U}}(\hbar\gamma)}}(\g_2[1])$.

C'est-\`a-dire pour tout couple $(\alpha,\beta)$ de poly-champs de
vecteurs tels que $[\alpha,\gamma]_{SN}=[\beta,\gamma]_{SN}=0$ on
a
\begin{equation}
d{\cal U}^{0}(\alpha\cup\beta)=d{\cal U}^{0}(\alpha)\cup d{\cal
U}^{0}(\beta)+D,
\end{equation}
o\`u $D$ est un cobord de Hochschild de l'alg\`ebre
$(C^{\infty}(X)[[\hbar]],\star_\gamma)$ donn\'e par $$
D=-[\star_\gamma,d{\cal U}^1(\alpha,\beta)]_G.
$$
\end{thm}\vskip 15pt
De fa\c con plus g\'en\'erale, l'\'el\'ement d'homotopie $D$ est
donn\'e par \begin{equation} D=-[\star_\gamma,d{\cal
U}^1(\alpha,\beta)]_G+d{\cal
U}^1(Q^{\hbar\gamma}(\alpha\otimes\beta)),
\end{equation}
avec pour $\alpha \in T^{|\alpha|}_{\rm{poly}}(X) $ et $\beta \in
T^{|\beta|}_{\rm{poly}}(X)$\begin{equation}\nonumber
Q^{\hbar\gamma}(\alpha\otimes\beta)=-[\hbar\gamma,
\alpha]_{SN}\otimes \beta-(-1)^{|\alpha|+1}\alpha\otimes
[\hbar\gamma,\beta]_{SN}.
\end{equation}
 \section{Isomorphisme de Duflo en cohomologie}
\subsection{Quantification du crochet de
Kirillov-Kos\-tant-Poisson}
 Dans le cas o\`u la vari\'et\'e $X$
est le dual d'une alg\`ebre de Lie de dimension finie $\g$ les
c\oe fficients du 2-champ de vecteurs de Kirillov-Kostant-Poisson
$\gamma$ sont des fonctions lin\'eaires sur $\g^*$. Si
$\{e_1,\dots,e_d\}$ est une base de $\g$ et $(e_1^*, \ldots,
e_d^*)$ sa base duale  on notera
$$\gamma=\frac12\sum_{i,j}[e_i,e_j]e_i^*\wedge e_j^*$$le 2-champ de vecteurs de Poisson associ\'e.

Ceci \'etant on peut consid\'erablement simplifier l'expression
des op\'erateurs poly-diff\'erentiels $B_{\Gamma}$ qui vont
intervenir dans la d\'efinition du star-produit $\star_\gamma$
(formule (\ref{starproduit})). A cause de la lin\'earit\'e des
c\oe fficients de $\gamma$ seuls interviennent les graphes dont
les sommets du premier type re\c coivent au plus une ar\^ete
(chaque sommet a\'erien est le but d'au plus une ar\^ete). On
parlera alors de graphes admissibles lin\'eaires.
Leur structure est bien d\'ecrite et relativement simple (voir par exemple \cite{[AST]}).\\

Il est facile de voir alors que si $f_1$ et $f_2$ sont deux
polyn\^omes alors $f_1\star_{\gamma}f_2$ est en fait une somme
finie. En localisant en $\hbar =1$ on obtient un  produit
associatif sur l'alg\`ebre des polyn\^omes sur $\g^*$. Ce produit
d\'efinit sur $S(\g)$ une structure d'alg\`ebre isomorphe \`a
l'alg\`ebre enveloppante $U(\g)$.\\

 Dans \cite{[K97]} \S  8.4,  Kontsevich
note $I_{alg}$ l'isomorphisme d'alg\`ebres entre $(S(\g),
\star_\gamma)$ et $U(\g)$. Compte tenu d'un r\'esultat de Shoikhet
\cite{[Sh3]} le r\'esultat 8.3.4 de \cite{[K97]} montre que
l'identification entre $(S(\g), \star_\gamma)$ et $U(\g)$ se fait
\`a l'aide de l'isomorphisme (d'espaces vectoriels)  de Duflo.\\

Rappelons ce qu'est la formule de Duflo. On consid\`ere   la
s\'erie formelle sur $\g^*$ d\'efinie au voisinage de $0$  pour
$x\in \g$ par
\begin{equation}q(x)=\det_{\g}\left(\frac{\sinh
(\frac{\mathrm{ad}x}{2})} {\frac{\mathrm{ad}x}{2}}\right)^{\frac
12}.\end{equation}On note $\partial(q)$ l'op\'erateur
diff\'erentiel (d'ordre infini) sur $S(\g)$ correspondant et on
note $\beta$ la sym\'etrisation de $S(\g)$ dans $U(\g)$.
L'isomorphisme (d'espaces vectoriels) de Duflo de $S(\g)$ dans
$U(\g)$ s'\'ecrit
\begin{equation}\label{duflo}\beta\circ
\partial(q).
\end{equation}Le point remarquable est que cette application restreinte aux invariants est un
isomorphisme d'alg\`ebres. On verra que c'est encore le cas pour
toute la cohomologie.
\subsection{Complexes de Chevalley-Eilenberg et de Hochschild}
En ne consid\'erant que les objets polynomiaux on peut sans
difficult\'es localiser en $\hbar=1$.

Le th\'eor\`eme (\ref{homotopie}) rappel\'e ci-dessus se restreint
alors au cadre poly\-nomial. On dispose donc d'un isomorphisme (la
diff\'erentielle $d\mathcal{U}^0$) entre l'espace de cohomologie
de Poisson de $S(\g)$ (c'est \`a dire la cohomologie du complexe
des poly-champs de vecteurs sur $\g^*$ \`a c\oe fficients dans
$S(\g)$ munie du cobord $[\gamma,-]_{SN}$ ) et la cohomologie de
Hochschild de l'alg\`ebre enveloppante $U(\g)$ munie du cup
produit standard. On cherche \`a expliciter cet isomorphisme.
\vskip 10pt
 Afin de simplifier les calculs
nous \'etablissons une autre r\'ealisation de ce dernier espace de
cohomologie.\\

En g\'en\'eral, \'etant donn\'e une alg\`ebre de Lie de dimension
finie $\g$ et un $\g$-module \`a gauche $M$ on d\'efinit la
cohomologie $H^*(\g,M)$ de $\g$ \`a c\oe fficients dans $M$ comme
le foncteur d\'eriv\'e \`a droite du foncteur  des invariants.
Plus concrètement ({\it cf}  \cite {[CE],[ChE],[We]}), les modules
de cohomologie $$H^*(\g,M)$$se calculent par  la cohomologie du
complexe des cocha\^\i nes de Chevalley-Eilenberg
$$\Hom_{\g}(V(\g),M),$$ o\`u $V(\g)$ d\'esigne le complexe
standard de Chevalley-Eilenberg. Rappelons que l'on a
$V_p(\g)=U(\g)\otimes\bigwedge^p\g$ avec diff\'erentielle
$$d:V_p(\g)\to V_{p-1}(\g)$$  donn\'ee par
\begin{eqnarray*}
&&\hspace{-0,5cm}d(u\otimes e_1\wedge\dots\wedge
e_p)=\sum_{i=1}^p(-1)^{i+1}ue_i\otimes e_1\wedge\dots\wedge
e_{i-1}\wedge e_{i+1}\wedge\dots\wedge e_p\\&&\hspace{-1cm}
+\sum_{i<j}(-1)^{i+j}u\otimes[e_i,e_j]\wedge e_1\wedge\dots\wedge
e_{i-1}\wedge e_{i+1}\wedge\dots\wedge e_{j-1}\wedge
e_{j+1}\wedge\dots\wedge e_p.
\end{eqnarray*}
Le complexe des cocha\^\i nes de Chevalley-Eilenberg
$$
\Hom_{\g}(V(\g),M)=\Hom_{\g}(U(\g)\otimes\bigwedge \g,M)
$$
est isomorphe au complexe des cocha\^\i nes
\begin{equation}\label{complex}
\Hom_{\Bbb R}(\bigwedge\g,M)=\bigwedge\g^*\otimes M
\end{equation}
 dont le cobord $\delta:\bigwedge^n\g^*\otimes
M\to\bigwedge^{n+1}\g^*\otimes M$ est d\'efini par
\begin{eqnarray*}
&&\delta f(e_1,\dots,e_{n+1})=\sum_{i}(-1)^{i+1}e_if(e_1,\dots,e_{i-1},e_{i+1},\dots,e_{n+1})+\\
&+&\sum_{i<j}(-1)^{i+j}f([e_i,e_j],\dots,e_{i-1},e_{i+1},\dots,e_{j-1},e_{j+1},\dots,e_{n+1}).
\end{eqnarray*}\\
D'un autre c\^ot\'e, \`a toute alg\`ebre $A$ et \`a tout
$A$-bimodule $M$ on associe un module cosimplicial
$[n]\mapsto\Hom(A^{\otimes n},M)$ en posant:

\begin{eqnarray*}
&&\partial^i(f)(a_0,\ldots,a_n)=
\begin{cases}a_0 f(a_1,\ldots,a_n)&\rm{
si}\;i=0
\\f(a_0,\dots,a_{i-1}a_i,\dots,a_n)&\rm{pour}\;0<i\leq
n
\\f(a_0,\dots,a_{n-1})a_n&\rm{si}\; i=n+1\end{cases}
\\
&&(\sigma^if)(a_1,\ldots,a_{n-1})=f(a_1,\dots,a_i,1,a_{i+1},\dots,a_{n-1}).
\end{eqnarray*}
La cohomologie de Hochschild $HH^*(A,M)$ de $A$ \`a c\oe fficients
dans $M$ ce sont les modules
$$
HH^n(A,M)=H^nC(\Hom(A^{\otimes },M)),
$$
o\`u $C\Hom(A^{\otimes },M)$ est le complexe des cocha\^\i nes
associ\'e:
$$
0\to
M\stackrel {\partial_0-\partial_1}\longrightarrow \Hom(A,M)\stackrel d\longrightarrow\Hom(A\otimes
A,M)\stackrel d\longrightarrow\dots
$$
dont le cobord $d$ est donn\'e par $d=\sum(-1)^i\partial^i$.
\\

Dans le cas o\`u on a $A=U(\g)$  et $M=U(\g)$ on obtient le
r\'esultat suivant ({\it cf} \cite{[Lo]} Lemma 3.3.3) :
\begin{lem}\label{lemma}
Les espaces de cohomologie $HH^*(U(\g),U(\g))$ et $H^*(\g,U(\g))$
sont isomorphes par l'anti-sym\'etriseur.
\end{lem}

\noindent Cet isomorphisme est induit par le morphisme de
complexes :
$$
\Psi^*\::\:\Hom_\R(U(\g)^{\otimes},
U(\g))\mapsto\Hom_\R(\bigwedge\g,U(\g))
$$
d\'efini pour $f\in \Hom(U(\g)^{\otimes n}, U(\g))$ par
\begin{eqnarray*}
(\Psi^* f)( e_1\wedge\dots\wedge e_n)&=&f(\Psi(
e_1\wedge\dots\wedge e_n))\\&:=&f\left(\frac{1}{n!}\sum_{\sigma\in
S_n}\sgn(\sigma) e_{\sigma(1)}\otimes\ldots\otimes
e_{\sigma(n)}\right).\end{eqnarray*} On a alors
\begin{eqnarray*}
\Psi^*([\star,f]_G)&=&\delta(\Psi^*(f)).
\end{eqnarray*}
\noindent
 Ainsi nous pouvons identifier le module de cohomologie de
Hochschild tangentielle de l'alg\`ebre de Lie diff\'erentielle
gradu\'ee $\g_2[1]$ avec celui du complexe (\ref{complex}). Cela
veut dire que pour comprendre la diff\'erentielle du
$L_\infty$-morphisme il faut faire agir les \'el\'ements de $
H_{\overline{\mathcal{U}}(\gamma)}$ sur les fonctions lin\'eaires
et prendre l'anti-sym\'etrisation.\\

En conclusion, \'etant donn\'e le morphisme d'alg\`ebres $d{\cal
U}^0$ d\'efini par (\ref{uderiv}) nous disposons d'un morphisme
d'alg\`ebres par composition :
$$
H^*_{\rm{Poisson}}(\g,S(\g))\stackrel {d{\cal U}^0}\longrightarrow
HH^*(U(\g),U(\g))\stackrel {\Psi^*}\longrightarrow H^*(\g,U(\g)).
$$Le reste de cette note est consacr\'e \`a la formule explicite pour ce
morphisme d'alg\`ebres.

\subsection{Les graphes intervenant dans  $d{\cal
U}^0$}

 Nous allons \`a pr\'esent calculer
explicitement l'expression du morphisme $d{\cal
U}^{0,\gamma}=d{\cal U}^0$ (formule (\ref{uderiv})) dans le cas
o\`u la vari\'et\'e $X$ est le dual d'une alg\`ebre de Lie
r\'eelle de dimension finie munie du $2-$champ de vecteurs de
Kirillov-Kostant-Poisson $\gamma$.\\

Les graphes qui vont intervenir dans l'expression du morphisme
$d{\cal U}^0$ sont $0-$admissibles.  \'Etant donn\'ee
l'identification du module de la cohomologie de Hochschild
tangentielle avec la cohomologie du complexe des cocha\^\i nes de
Chevalley-Eilenberg (Lemme \ref{lemma}),  on doit \'evaluer
$d{\cal U}^0(\alpha)$ avec $\alpha\in T_{\rm{poly}}(\g^*)[1]$ sur
des fonctions lin\'eaires.\\

\begin{figure}[!h]
\begin{center}
  % Requires \usepackage{graphicx}
 \includegraphics[height=1.7in]{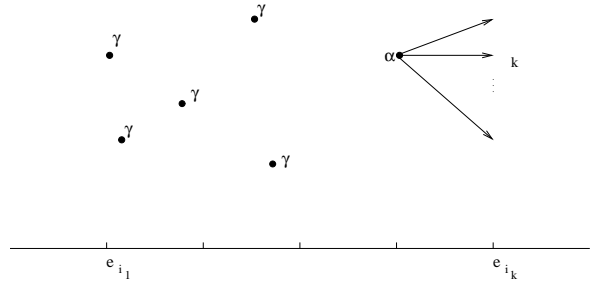}\\
\caption{Les graphes de $d{\cal U}^0(\alpha)$ }\label{2}
\end{center}
\end{figure}

\begin{figure}[!h]
\begin{center}
  % Requires \usepackage{graphicx}
 \includegraphics[width=2.7928in]{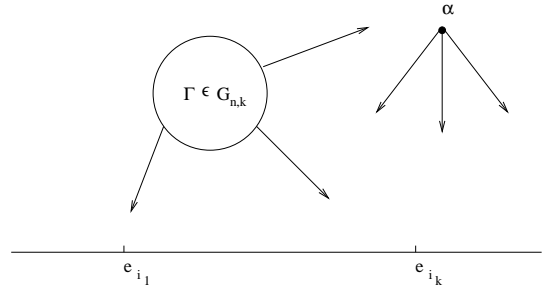}\\
  \caption{Les graphes de $d{\cal U}^0(\alpha)$ }\label{3}
  \end{center}
\end{figure}

Les  graphes et les op\'erateurs qui vont intervenir sont comme
sur les figures [\ref{2}, \ref{3}] o\`u l'on associe $k$ fonctions
lin\'eaires de coordonn\'ees $e_{i_1}, \ldots e_{i_k}$ aux points
terrestres, le $k-$champ de vecteur $\alpha$ au point de
$V_{\Gamma}^{1,1}=\{1\}$ et le $2-$champ de vecteurs  $\gamma$ aux
$n$ points a\'eriens restants.

Pour que le poids correspondant ne s'annule pas, il faut qu'il y
ait exactement $2n+k$ ar\^etes, ce qui correspond bien \`a la
d\'efinition des graphes $0$-admissibles.\\

En utilisant la lin\'earit\'e des c\oe fficients du 2-champ de
vecteurs $\gamma$, nous allons d\'eterminer la forme exacte de
tels graphes.\\

Soit $\widehat{\Gamma}$ un graphe contribuant a priori dans le
calcul de $d{\cal U}^0(\alpha)$. D\'esignons par $\Gamma$ le sous
graphe dans lequel on a enlev\'e le sommet $\alpha$ et les $k$
ar\^etes issues de ce sommet (figure [\ref{3}]).

Le graphe $\Gamma$ se d\'ecompose en r\'eunion de graphes simples
qui sont soit de type Lie (c'est \`a dire
des arbres) soit de type roue tentaculaire.\\

 Les ar\^etes qui
partent du sommet de $V_{\Gamma}^{1,1}$ que l'on appellera sommet
$\alpha$ vont donc soit sur les sommets terrestres soit sur les
racines des sous-arbres ind\'ependants du graphe $\Gamma$.

\'Etant donn\'e que les fonctions $e_{i_j}$ "d\'eriv\'ees" par le
graphe $\Gamma$ sont lin\'eaires, le seul point de "connexion" des
sous-arbres ind\'ependants est donc le sommet $\alpha$ :
\begin{figure}[!h]
\begin{center}

 \includegraphics[width=3in]{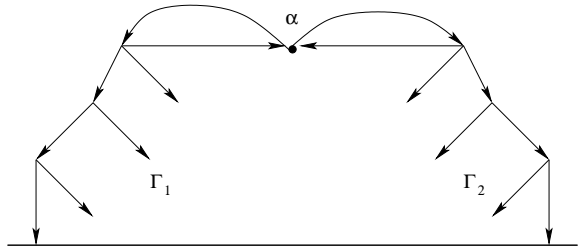}\\
  \caption{Point de connexion}\label{pointconnexion}
  \end{center}
\end{figure}

Ceci montre que la  forme $\Omega_{\widehat{\Gamma}}$ se
pr\'esente comme  produit s\'epar\'e de formes associ\'ees aux
sous-graphes simples
 de $\Gamma$ auxquels on a ajout\'e \'eventuellement une ar\^ete issue de $\alpha$.
On est donc amen\'e \`a \'etudier les contributions dans
les formes d'angles des graphes simples de $\Gamma$ \'etendus \'eventuellement par une ar\^ete issue de $\alpha$.\\

\begin{figure}[ptbh]
\begin{center}
  % Requires \usepackage{graphicx}
 \includegraphics[width=1.1in]{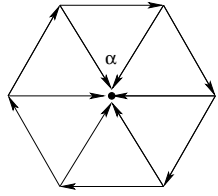}\\
  \caption{Une roue pure}\label{roue}
  \end{center}
\end{figure}

\noindent {\bf Premier cas :} Soit $\Gamma_1$ un sous-graphe
simple de $\Gamma$ ne recevant pas de fl\`eches provenant du
sommet $\alpha$. Supposons que ce graphe poss\`ede $p$ sommets
a\'eriens et $k$ sommets terrestres. Ce graphe poss\`ede $2p$
ar\^etes. La forme d'angle associ\'ee est int\'egr\'ee sur une
vari\'et\'e de configurations  de dimension $2p+k$ on en d\'eduit
que $k=0$. Alors au plus $p$ fl\`eches d\'erivent les sommets
a\'eriens (correspondant au $2$-vecteur lin\'eaire $\gamma$) et au
moins $p$ fl\`eches issues des sommets a\'eriens d\'erivent le
poly-champ de vecteurs $\alpha$. Or, les ar\^etes doubles et les
boucles \'etant interdites, de chacun de $p$ sommets a\'eriens
part une et une seule ar\^ete vers  le sommet $\alpha$. Ainsi on
se retrouve n\'ecessairement avec un graphe de type roue pure
comme \`a la figure [\ref{roue}]. La contribution dans le poids de
ce sous-graphe se factorise d'une part et d'autre part est nulle
d'apr\`es \cite{[Sh3]}.

\noindent{\bf Deuxi\`eme cas: } Consid\'erons $\Gamma_1$ un
sous-graphe simple de racine num\'erot\'e $1$ et recevant une
ar\^ete $e=(\alpha, 1)$ du sommet $\alpha$. Notons $p_1$ le nombre
de ses points a\'eriens et $k_1$ le nombre de ses points
terrestres. Le graphe \'etendu $\widehat{\Gamma_1}=\Gamma_1\cup
\{e\}$ a donc $2p_1+1$ ar\^etes. Or la contribution  dans le
calcul total du poids, s'obtient par intégration de la forme
$\Omega_{\widehat{\Gamma_1}}$ sur un espace de configurations de
dimension $2p_1+k_1$.

On a donc $k_1=1$ et le sommet terrestre est effectivement atteint
sinon  la dimension de l'espace de configurations serait
inf\'erieure \`a $2p_1+1$.

Le graphe $\Gamma_1$ poss\`ede alors $2p_1$ ar\^etes dont au plus
$p_1-1$ vont sur les points a\'eriens de $\Gamma_1$ (c'est un
arbre) et au moins $p_1$ autres vont  sur $\alpha$  et une arête
va sur un sommet terrestre. Comme les ar\^etes doubles sont
interdites, on en d\'eduit que de chacun des $p_1$ sommets de
$\Gamma_1$ part une et une seule ar\^ete vers le sommet $\alpha$.
Il est alors presque \'evident que l'arbre $\Gamma_1$ est de la
forme de la figure [\ref{escargot}]. On appellera une telle figure
un escargot.\\

En reprenant les notations introduites plus haut on peut donc
\'enoncer le lemme suivant :

\begin{lem} Soit  $\widehat{\Gamma}$ un graphe 0-admissible dont la contribution
dans la formule de la formalit\'e tangente est a priori non nulle.
Alors le graphe $\Gamma$ se d\'ecompose est un  produit de graphes
simples comme \`a la figure [\ref{escargot}], appel\'es escargots.
\begin{figure}[ptbh]
\begin{center}
  % Requires \usepackage{graphicx}
\includegraphics[width=3in]{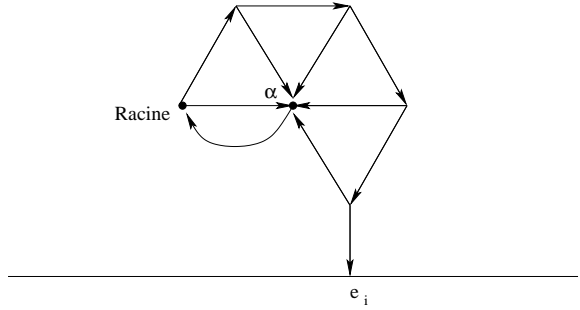}\\
  \caption{Un escargot}\label{escargot}
  \end{center}
\end{figure}
\end{lem}

\subsection{Factorisation des contributions des escargots}
 \'Etant donn\'e un $0$-graphe admissible quelconque sa forme d'angle associ\'ee
 est donc le produit des formes d'angles  des escargots qui le
composent.  Nous allons voir que l'antisym\'etrisation \'evoqu\'ee
dans le lemme \ref{lemma} permet de factoriser dans le poids les
contributions de ces graphes
ce qui nous permet d'\'enoncer le lemme suivant.\\

\begin{lem}
Les contributions dans le poids des produits d'escargots se
factorisent gr\^ace \`a l'op\'eration d'antisym\'etrisation.
\end{lem}

\begin{figure}[!h]
\begin{center}

 \includegraphics[width=2.2in]{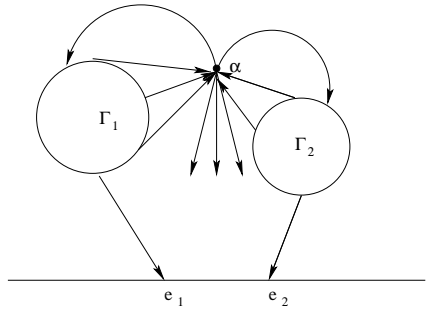}\\
  \caption{Produit de deux escargots}\label{deuxescargots}
  \end{center}
\end{figure}

Raisonnons dans le cas o\`u $\alpha$ est un $p-$champ de vecteurs.
Les graphes qui interviennent sont des produits de $p$ escargots
(\'eventuel\-lement triviaux), not\'e $\widehat{\Gamma_1}, \ldots,
\widehat{\Gamma_p}$ comme \`a la figure [\ref{deuxescargots}] o\`u
on n'a repr\'esent\'e que deux escargots. Ces graphes sont
rep\'er\'es par la position terrestre correspondante et sont
rang\'es dans l'ordre croissant.

\begin{figure}[!h]
\begin{center}

 \includegraphics[width=2.2in]{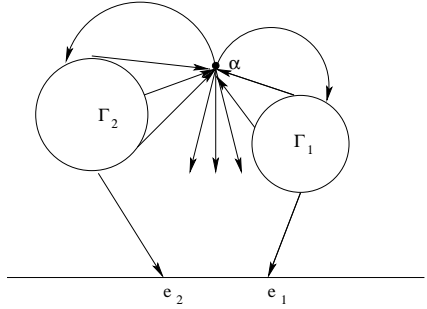}\\
  \caption{L'anti-sym\'etris\'e du graphe associé à $\widehat{\Gamma}$}\label{antisym}
  \end{center}
\end{figure}

Appliquons  l'antisym\'etrisation qui nous permet d'identifier les
\'el\'e\-ments de la cohomologie tangentielle de Hochschild avec
ceux de la cohomologie des cocha\^\i nes sur $U(\g)$ (Lemme
\ref{lemma}). Nous obtenons des contributions de la forme
$$\sum_{\sigma} \frac{\sgn(\sigma)}{p!} \left(\int \Omega_{\widehat{\Gamma_1}}
\wedge\ldots\wedge \Omega_{\widehat{\Gamma_p}}\right) \;
B_{\widehat{\Gamma_1}\cup\ldots\cup\widehat{\Gamma_p}}(\alpha)(e_{\sigma(1)}
 \otimes\ldots\otimes e_{\sigma(p)}).$$
Consid\'erons maintenant les contributions des graphes
sym\'etris\'es de
$\widehat{\Gamma_1}\cup\ldots\cup\widehat{\Gamma_p}$. Par graphe
sym\'etris\'e on entend le graphe obtenu en  permutant les
positions terrestres  des graphes $\widehat{\Gamma_i}$ comme à la
figure [\ref{antisym}]. On obtient pour toute permutation $\tau$
les contributions de la forme
\begin{equation}\label{contribution}\hspace{-0,5cm}\sum_{\sigma}\frac{\sgn(\sigma)}{p!}
\left(\int \Omega_{\widehat{\Gamma_{\tau(1)}}} \wedge\ldots\wedge
\Omega_{\widehat{\Gamma_{\tau(p)}}}\right) \;
B_{\widehat{\Gamma_{\tau(1)}}\cup\ldots\cup\widehat{\Gamma_{\tau(p)}}}(\alpha)(e_{\sigma(1)}
 \otimes\ldots\otimes e_{\sigma(p)}).
 \end{equation}

Clairement on a
\begin{equation}\label{symetrisation}B_{\widehat{\Gamma_{\tau(1)}}\cup\ldots\cup
\widehat{\Gamma_{\tau(p)}}}(\alpha)(e_{\tau(1)}\otimes\ldots\otimes
e_{\tau(p)})=B_{\widehat{\Gamma_1}\cup\ldots\cup\widehat{\Gamma_p}}(\alpha)(e_1\otimes\ldots\otimes
e_p).
\end{equation}\\

Notons $x_i$ la position terrestre correspondant au graphe
$\widehat{\Gamma}_i$. La contribution dans le poids de l'escargot
$\widehat{\Gamma_i}$ obtenue par int\'egration sur les positions
de
 premi\`ere
 esp\`ece (points a\'eriens) o\`u est plac\'e le $2$-vecteur de Poisson est donc une 1-forme en
 $x_i$, not\'ee $\psi_{\widehat{\Gamma_i}}(x_i)dx_i$. On a donc :

 \begin{equation}\label{forme}
\left(\int \Omega_{\widehat{\Gamma_{\tau(1)}}} \wedge\ldots\wedge
\Omega_{\widehat{\Gamma_{\tau(p)}}}\right)= \int_{x_1<\ldots<
x_p}\psi_{\widehat{\Gamma_1}}(x_1)\ldots
\psi_{\widehat{\Gamma_p}}(x_p)dx_1\wedge\ldots\wedge
dx_p.\end{equation}

 En sommant les diff\'erentes contributions (\ref{contribution})
pour chaque permutation $\tau$ on obtient, compte tenu des
\'equations (\ref{symetrisation}) et (\ref{forme}), pour tout
$\widehat{\Gamma}_1, \ldots, \widehat{\Gamma}_p$ : $$\sum_{\tau,
\sigma}\frac{\sgn(\sigma)}{p!} \left(\int
\Omega_{\widehat{\Gamma_{\tau(1)}}} \wedge\ldots\wedge
\Omega_{\widehat{\Gamma_{\tau(p)}}}\right)\cdot
B_{\widehat{\Gamma_{\tau(1)}}\cup\ldots\cup\widehat{\Gamma_{\tau(p)}}}(\alpha)(e_{\sigma(1)}
 \otimes\ldots\otimes e_{\sigma(p)})=$$
 \begin{eqnarray}\label{final}
 &&\hspace{-1cm}\sum_{\tau, \sigma}\frac{\sgn(\sigma)\sgn(\tau)}{p!}
\left(\int_{x_{\tau(1)}<\ldots<x_{\tau(p)}}
\psi_{\widehat{\Gamma_{\tau(1)}}}(x_{\tau(1)})dx_{\tau(1)}\wedge
\ldots
\wedge\psi_{\widehat{\Gamma_{\tau(p)}}}(x_{\tau(p)})dx_{\tau(p)}\right)\cdot\nonumber\\
&&\hspace{5cm}B_{\widehat{\Gamma_1}\cup\ldots\cup\widehat{\Gamma_p}}(\alpha)(e_{\sigma(1)}
 \otimes\ldots\otimes e_{\sigma(p)}).
\end{eqnarray}

  Or on a :
 \begin{eqnarray}&&\hspace{-1cm}\sum_{ \tau}\sgn(\tau)
\left(\int_{x_{\tau(1)}<\ldots<x_{\tau(p)}}
\psi_{\widehat{\Gamma_{\tau(1)}}}(x_{\tau(1)})dx_{\tau(1)}\wedge
\ldots
\wedge\psi_{\widehat{\Gamma_{\tau(p)}}}(x_{\tau(p)})dx_{\tau(p)}\right)=\nonumber\\&&\hspace{3cm}
\int_{x_1, \ldots, x_p}\psi_{\widehat{\Gamma_1}}(x_1)\ldots
\psi_{\widehat{\Gamma_p}}(x_p)dx_1\wedge\ldots\wedge
dx_p,\end{eqnarray} qui est une int\'egrale \`a variables
s\'epar\'ees.  Au final l'expression (\ref{final}) s'\'ecrit
\begin{equation}\label{somme}
\sum_{
\sigma}\frac{\sgn(\sigma)}{p!}\left(\prod_{i}\int_{\R}\psi_{\widehat{\Gamma_i}}(x_i)dx_i\right)B_{\widehat{\Gamma_1}\cup\ldots\cup\widehat{\Gamma_p}}(\alpha)(e_{\sigma(1)}
 \otimes\ldots\otimes e_{\sigma(p)}).
 \end{equation}

\subsection{Calcul du poids d'un escargot}

\begin{lem}\label{annulation}
Le poids associ\'e \`a un escargot non trivial est nul. \end{lem}

Soit $\widehat{\Gamma}=\Gamma\sqcup (1,\alpha)$ un escargot ({\it
i.e} un arbre 0-admissible) avec $p+1$ sommets a\'eriens,
c'est \`a dire que l'arbre $\Gamma$ a $p$ sommets a\'eriens.\\

Fixons la position du sommet terrestre en $0$ tandis que  celle du
sommet $\alpha$ est rep\'er\'e sur le demi-cercle unit\'e  par
l'angle $\theta$. L'espace de configurations est fibr\'e par la
position du sommet $\alpha$, et chaque fibre est une
sous-vari\'et\'e de dimension $2p$. Notons
$w_{\Gamma}(\theta)d\theta$ l'int\'egrale de la forme
diff\'erentielle $\Omega_{\widehat{\Gamma}}$ sur cette fibre.
C'est une 1-forme en la \mbox{variable $\theta$}.
\begin{figure}[ptbh]
\begin{center}
  % Requires \usepackage{graphicx}
 \includegraphics[width=2.2in]{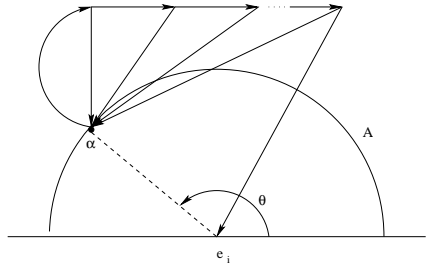}\\
  \caption{Calcul du poids}\label{calcul}
  \end{center}
\end{figure}

\begin{figure}[!h]
\begin{center}
  % Requires \usepackage{graphicx}
 \includegraphics[width=2.2in]{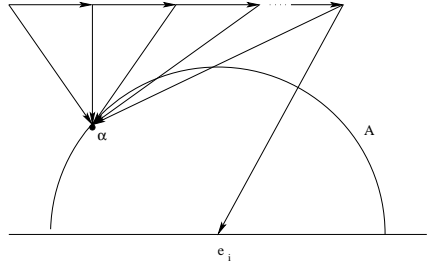}\\
  \caption{Graphe d\'eploy\'e }\label{deploye}
\end{center}
\end{figure}

 Le poids de l'escargot
$\widehat{\Gamma}$ est donc l'int\'egrale
$w_{\Gamma}(\theta)d\theta$ le long de l'arc de demi-cercle:
$$w_{\widehat{\Gamma}}=\int_0^{\pi}w_{\Gamma}(\theta)d\theta.$$Or la $1-$forme
diff\'erentielle $w_{\Gamma}(\theta)d\theta$ est la
diff\'erentielle du poids du graphe d\'eploy\'e
$\overset{\vee}{\Gamma}$ (on a d\'eploy\'e l'ar\^ete issue du
sommet $\alpha$ en ajoutant un sommet). C'est un graphe avec $p+2$
points a\'eriens et $2p+2$ ar\^etes comme \`a la figure
[\ref{deploye}].

En effet, consid\'erons la sous-vari\'et\'e des configurations \`a
$\alpha$ fixe. C'est une sous-vari\'et\'e de dimension $2p+2$. On
note alors $w_{\overset{\vee}{\Gamma}}(\theta)$ l'int\'egrale de
la $2p+2$ forme $\Omega_{\overset{\vee}{\Gamma}}$ sur cette
sous-vari\'et\'e.

Le calcul de la diff\'erentielle de
$w_{\overset{\vee}{\Gamma}}(\theta)$ se fait par contraction
d'ar\^etes. C'est une cons\'equence de la formule de Stokes  et du
lemme 6.6 de \cite{[K97]} (dont on peut se passer dans notre cas,
car les graphes sont lin\'eaires). On voit alors facilement que
cette diff\'erentielle est pr\'ecis\'ement
$\pm w_{\Gamma}(\theta)d\theta$ (voir \cite{[T]}).\\

Le graphe $\overset{\vee}{\Gamma}$ est alors un graphe de type
Bernoulli ({\it cf} \cite {[Kat], [T]}) avec $p+2$ sommets
a\'eriens. Les c\oe fficients $w_{\overset{\vee}{\Gamma}}(0)$ et
$w_{\overset{\vee}{\Gamma}}(\pi)$ correspondent donc au poids des
graphes de Bernoulli:
\begin{figure}[!h]
\begin{center}
  % Requires \usepackage{graphicx}
\includegraphics[height=1.1in,width=3.5in]{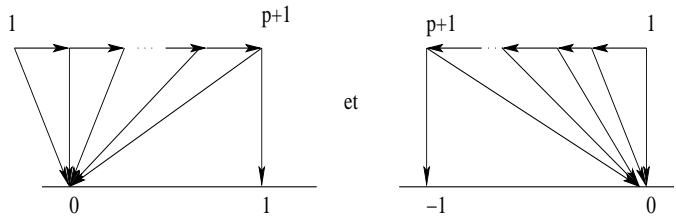}\\
  \caption{Graphes de Bernoulli}\label{bernoulli}
\end{center}
\end{figure}

\noindent qui valent respectivement
$\frac{(-1)^{p+1}}{(p+1)!}b_{p+1}(0)$ et
$\frac{(-1)^{p+1}}{(p+1)!}b_{p+1}(1)$ o\`u $b_n(x)$ d\'esigne le
$n$-i\`eme polyn\^ome de Bernoulli ({\it cf} \cite {[Kat]}
\S4.4.3, \cite{[T]} \S4). Or, d'un c\^ot\'e les polyn\^omes de
Bernoulli v\'erifient l'\'egalit\'e
$$
b_n(1)=(-1)^nb_n(0)
$$
et de l'autre on sait que tous les nombres de Bernoulli
$B_n=b_n(0)$ sont nuls pour $n$ impair plus grand que $3$. Pour
$p>0$ on  conclut que l'on a
$$w_{\widehat{\Gamma}}=w_{\overset{\vee}{\Gamma}}(0)-w_{\overset{\vee}{\Gamma}}(\pi)=0,$$et
pour $p=0$ on retrouve bien s\^ur la valeur $1=\frac 12 -(-\frac
12)$.

\subsection{Conclusion}
 Nous pouvons donc conclure que l'application  $\Psi^* \circ
 d\mathcal{U}^0$
 est
 en fait l'identit\'e, lorsque l'on voit $U(\g)$ comme l'alg\`ebre $(S(\g), \star)$ via l'application de Duflo
 (formule (\ref{duflo})).\\

En effet, d'apr\`es le lemme \ref{annulation} seuls les escargots
triviaux  vont
 intervenir dans la formule (\ref{somme}) et le c\oe fficient
 $$\prod_{i}\int_{\R}\psi_{\widehat{\Gamma_i}}(x_i)dx_i$$ vaut alors clairement $1$.
 Pour  tout  $p-$champ de vecteurs \`a c\oe fficients
 polynomiaux $\alpha$ on obtient donc
 \begin{eqnarray}
 \Psi^*\big(
 d\mathcal{U}^0(\alpha)\big)(e_1\wedge\ldots\wedge e_p)&=&\frac{1}{p!}\sum_{\sigma}\sgn(\sigma)\alpha(e_{\sigma(1)}
  \otimes\ldots\otimes e_{\sigma(p)})\nonumber\\&=&\alpha(e_1\wedge\ldots\wedge
  e_p).
  \end{eqnarray}
  On peut énoncer le résultat de cette note :

\begin{thm}
L'application de Duflo s'\'etend en un isomorphisme d'alg\`ebres
de $H_{Poisson}(\g,S(\g))$ sur $H(\g,U(\g)).$
\end{thm}
%\newpage
\footnotesize{}
\end{document}